\documentclass[submission]{mydmtcs}

\usepackage[latin1]{inputenc}
\usepackage{subfigure}

\usepackage[centertags]{amsmath}
\usepackage{amsfonts}
\usepackage{amssymb}
\usepackage{newlfont}
\usepackage{color}

\author{Michael Huber}

\title[Block-transitive combinatorial designs]
{On the existence of block-transitive combinatorial designs}

\address{Wilhelm-Schickard-Institute for Computer Science, University of Tuebingen, Sand~13,
D-72076 Tuebingen, Germany\\
E-mail \email{michael.huber@uni-tuebingen.de}}

\keywords{Combinatorial design, block-transitive group of automorphisms, \mbox{$3$-homogeneous} permutation group}

\received{21 Nov 2008}
\revised{18 Jan 2010}
\accepted{3 Mar 2010}

\mathchardef\ordinarycolon\mathcode`\:
  \mathcode`\:=\string"8000
  \begingroup \catcode`\:=\active
    \gdef:{\mathrel{\mathop\ordinarycolon}}
  \endgroup
\newtheorem{mthm}{Main Theorem}
\newtheorem{thm}{Theorem}

\newtheorem{prop}[thm]{Proposition}
\newtheorem{rem}[thm]{Remark}
\newtheorem{Cor}[thm]{Corollary}
\newtheorem{problem}[thm]{Problem}
\newtheorem{conj}[thm]{Conjecture}

\begin{document}

\maketitle

\begin{abstract}
Block-transitive Steiner \mbox{$t$-designs} form a central part of
the study of highly symmetric combinatorial configurations at the
interface of several disciplines, including group theory, geometry,
combinatorics, coding and information theory, and cryptography. The
main result of the paper settles an important open question:
There exist no non-trivial examples with $t=7$ (or larger).
The proof is based on the classification of the finite
\mbox{$3$-homogeneous} permutation groups, itself relying on the
finite simple group classification.
\end{abstract}

%%% ----------------------------------------------------------------------

\section{Introduction}\label{intro}

One of the outstanding problems in combinatorial design theory
concerns the existence of Steiner \mbox{$t$-designs} (\textit{i.e.},
\mbox{$t$-$(v,k,1)$} designs) with $t>5$. In particular the
existence of Steiner \mbox{$t$-designs} admitting an interesting
group of automorphisms is of great interest. The known examples for
$t \leq 5$ often encompass a high degree of regularity and establish
deep connections to permutation group theory, geometry, combinatorics,
coding and information theory, and cryptography.

There has been recent progress on the existence problem by
characterizing Steiner \mbox{$t$-designs} which admit a
flag-transitive group of automorphisms (\textit{cf.}~\cite{Hu2008}).
In this paper, we focus on the existence problem of Steiner
\mbox{$t$-designs} under the weaker condition of block-transitivity.
P.~Cameron and C.~Praeger~\cite{CamPrae1993} proved the
non-existence of block-transitive (Steiner) \mbox{$t$-designs} for $t > 7$.
Moreover, they conjectured that there are no non-trivial examples
for $t=6$. Recently, the author~\cite{Hu_mmics2008,Hu_cam2008}
essentially confirmed the non-existence of block-transitive Steiner
\mbox{$6$-designs}. The main result of the present paper settles now
the challenging remaining  question for $t=7$. We prove

\begin{mthm}\label{block7-des}
There is no non-trivial Steiner \mbox{$7$-design} $\mathcal{D}$ admitting a
block-transitive group \mbox{$G \leq \emph{Aut}(\mathcal{D})$} of automorphisms.
\end{mthm}

The paper is organized as follows: Preliminary results which are important for
the remainder of the paper are collected in Section~\ref{Prelim}.
In Section~\ref{prev}, a detailed account on previous and related work is presented.
In Section~\ref{proof}, the proof of the Main Theorem will be given. It is based on the classification
of the finite \mbox{$3$-homogeneous} permutation groups, itself
relying on the classification of the finite simple groups.

\smallskip

%%% ----------------------------------------------------------------------

\section{Preliminary Results}\label{Prelim}

\subsection{$t$-Designs}\label{designs}

Combinatorial design theory is a fascinating subject on the
interface of several disciplines, including group theory, geometry,
combinatorics, coding and information theory, and cryptography.
In particular, the study of designs with high symmetry properties has a
very long history and establishes deep connections between these
areas. One of its highlights surely is the remarkable interrelation
between the Mathieu--Witt designs, Golay codes, sporadic simple
Mathieu groups, Leech lattice, and Kissing Numbers and Sphere
Packing problems (\textit{cf.},~\textit{e.g.},~\cite{cosl98,Hu_cod2008,wisl77}).
A recent connection with cryptography is given in~\cite{Hub2009}.

Combinatorial designs may be regarded as generalizations of finite
projective planes. More formally: For positive integers $t \leq k
\leq v$ and $\lambda$, we define a \mbox{\emph{$t$-$(v,k,\lambda)$
design}} to be a finite incidence structure
\mbox{$\mathcal{D}=(X,\mathcal{B},I)$}, where $X$ denotes a set of
\emph{points}, $\left| X \right| =v$, and $\mathcal{B}$ a set of
\emph{blocks}, $\left| \mathcal{B} \right| =b$, satisfying the
following regularity properties: each block $B \in \mathcal{B}$ is
incident with $k$ points, and each \mbox{$t$-subset} of $X$ is
incident with $\lambda$ blocks. A \emph{flag} of $\mathcal{D}$ is an
incident point-block pair $(x,B) \in I$ with $x \in X$ and $B \in
\mathcal{B}$.

For historical reasons, a \mbox{$t$-$(v,k,\lambda)$ design} with
$\lambda =1$ is called a \emph{Steiner \mbox{$t$-design}} (sometimes
also a \emph{Steiner system}). We note that in this case each block
is determined by the set of points which are incident with it, and
thus can be identified with a \mbox{$k$-subset} of $X$ in a unique
way. If $t<k<v$, then we speak of a \emph{non-trivial} Steiner
\mbox{$t$-design}. As a simple example, the vector space
$\integers_2^n$ ($n \geq 3$) with block set $\mathcal{B}$ taken to
be the set of all subsets of four distinct elements of
$\integers_2^n$ whose vector sum is zero is a (boolean) Steiner
\mbox{$3$-$(2^n,4,1)$ design}. There are many infinite classes of
Steiner \mbox{$t$-designs} for $t=2$ and $3$, however for $t=4$ and
$5$ only a finite number are known. For a detailed treatment of
combinatorial designs, we refer
to~\cite{BJL1999,crc06,hupi85,stin04}. In
particular,~\cite{BJL1999,crc06} provide encyclopedic accounts of
key results and contain existence tables with known parameter sets.

We consider automorphisms of a \mbox{$t$-design} $\mathcal{D}$ as
pairs of permutations on $X$ and $\mathcal{B}$ which preserve
incidence, and call a group \mbox{$G \leq \mbox{Aut} (\mathcal{D})$}
of automorphisms of $\mathcal{D}$ \emph{block-transitive}
(respectively \emph{flag-transitive}, \emph{point
\mbox{$t$-transitive}}, \emph{point \mbox{$t$-homogeneous}}) if $G$
acts transitively on the blocks (respectively transitively on the
flags, \mbox{$t$-transitively} on the points,
\mbox{$t$-homogeneously} on the points) of $\mathcal{D}$. For short,
$\mathcal{D}$ is said to be, \textit{e.g.}, block-transitive if $\mathcal{D}$
admits a block-transitive group of automorphisms.

For \mbox{$\mathcal{D}=(X,\mathcal{B},I)$} a Steiner
\mbox{$t$-design} with \mbox{$G \leq \mbox{Aut} (\mathcal{D})$}, let
$G_x$ denote the stabilizer of a point $x \in X$, and $G_B$ the
setwise stabilizer of a block $B \in \mathcal{B}$. For $x, y \in X$
and $B \in \mathcal{B}$, we define $G_{xy}= G_x \cap G_y$.

For any $x \in \reals$, let $\lfloor x \rfloor$ denote the
greatest positive integer which is at most $x$.

All other notation is standard.

%%%------------------------------------------------------------------

\subsection{Combinatorial Existence Results}

We recall some standard combinatorial results which we use in this
paper. For the existence of \mbox{$t$-designs}, basic necessary
conditions can be obtained via elementary counting arguments (see,
for instance,~\cite{BJL1999}):

\begin{prop}\label{s-design}
Let $\mathcal{D}=(X,\mathcal{B},I)$ be a \mbox{$t$-$(v,k,\lambda)$}
design, and for a positive integer $s \leq t$, let $S \subseteq X$
with $\left|S\right|=s$. Then the total number of blocks incident
with each element of $S$ is given by
\[\lambda_s = \lambda \frac{{v-s \choose t-s}}{{k-s \choose t-s}}.\]
In particular, for $t\geq 2$, a \mbox{$t$-$(v,k,\lambda)$} design is
also an \mbox{$s$-$(v,k,\lambda_s)$} design.
\end{prop}
\noindent It is customary to set $r:= \lambda_1$ denoting the total
number of blocks incident with a given point (referring to the
`replication number' from statistical design of experiments, one of
the origins of design theory).

\begin{Cor}\label{Comb_t=5}
Let $\mathcal{D}=(X,\mathcal{B},I)$ be a \mbox{$t$-$(v,k,\lambda)$}
design. Then the following holds:
\begin{enumerate}

\item[{\emph(a)}] $bk = vr.$

\smallskip

\item[{\emph(b)}] $\displaystyle{{v \choose t} \lambda = b {k \choose t}.}$

\smallskip

\item[{\emph(c)}] $r(k-1)=\lambda_2(v-1)$ for $t \geq 2$.

\end{enumerate}
\end{Cor}

\begin{Cor}\label{divCond}
Let $\mathcal{D}=(X,\mathcal{B},I)$ be a \mbox{$t$-$(v,k,\lambda)$}
design. Then
\[\lambda {v-s \choose t-s} \equiv \, 0\; \emph{\bigg(mod}\;\, {k-s \choose t-s}\bigg)\]
for each positive integer $s \leq t$.
\end{Cor}

For non-trivial Steiner \mbox{$t$-designs} lower bounds for $v$ in
terms of $k$ and $t$ can be given (see
P.~Cameron~\cite[Thm.\,3A.4]{Cam1976}, and
J.~Tits~\cite[Prop.\,2.2]{Tits1964}):

\begin{prop}\label{Cam}
If $\mathcal{D}=(X,\mathcal{B},I)$ is a non-trivial Steiner
\mbox{$t$-design}, then the following holds:
\begin{enumerate}

\item[{\emph(a)}] \emph{(Tits~1964):} \hspace{0.2cm} $v\geq (t+1)(k-t+1).$

\smallskip

\item[{\emph(b)}] \emph{(Cameron~1976):} \hspace{0.2cm} \mbox{$v-t+1 \geq (k-t+2)(k-t+1)$} for $t>2$. If equality
holds, then
\smallskip
$(t,k,v)=(3,4,8),(3,6,22),(3,12,112),(4,7,23)$, or $(5,8,24)$.
\end{enumerate}
\end{prop}

\noindent We note that in general Part~(a) is stronger for
$k<2(t-1)$, while Part~(b) is stronger for \mbox{$k>2(t-1)$}. For
$k=2(t-1)$ both assert that $v \geq t^2-1$. As we are in particular
interested in the case when $t=7$, we deduce from Part~(b) the
following upper bound for the positive integer $k$.

\begin{Cor}\label{Cameron_t=5}
Let $\mathcal{D}=(X,\mathcal{B},I)$ be a non-trivial Steiner \mbox{$7$-design}. Then
\[k \leq \bigl\lfloor \sqrt{v} + \textstyle{\frac{11}{2}} \bigr\rfloor.\]
\end{Cor}

We finally state two classical results on the existence of
\mbox{$t$-designs}. The first is due to D.~Ray-Chaudhuri and
R.~Wilson~\cite[Thm.\,1]{Ray-ChWil1975}, and the second is by
L.~Teirlinck~\cite{Teir1987}:

\begin{thm}{\em (Ray-Chaudhuri \& Wilson~1975).}\label{RayCh}
Let $\mathcal{D}=(X,\mathcal{B},I)$ be a \mbox{$t$-$(v,k,\lambda)$}
design. If $t$ is even, say $t=2s$, and $v \geq k+s$, then $b \geq
{v\choose s}$. If $t$ is odd, say $t = 2s+1$, and $v-1 \geq k+s$,
then $b \geq 2{v-1\choose s}$.
\end{thm}

\begin{thm}{\em (Teirlinck~1987).}\label{Teirl}
For every positive integer value of $t$, there exists a non-trivial
\mbox{$t$-design}.
\end{thm}
\noindent However, although Teirlinck's recursive methods are
constructive, they only produce examples with tremendously large
values of $\lambda$. Via computer search, \mbox{$t$-$(v,k,\lambda)$}
designs with $t \geq 6$ and smaller values of $\lambda$ (where $\lambda$
is at least $4$) have been constructed in recent years by the method of
orbiting under a group (see~\cite{LaueKosh2007} for an overview).
Until now no non-trivial Steiner \mbox{$t$-design} with $t \geq 6$
is known.

\begin{problem}
Does there exist any non-trivial Steiner \mbox{$t$-design} with $t
\geq 6$?
\end{problem}

\smallskip

%%% ----------------------------------------------------------------------

\section{Previous and Related Work}\label{prev}

We focus on \mbox{$t$-designs} which admit groups of automorphisms
with high symmetry properties. One of the early important results is
due to R.~Block~\cite[Thm.\,2]{Block1965}:

\begin{prop}{\em (Block~1965).}\label{BlocksLemma2}
Let $\mathcal{D}=(X,\mathcal{B},I)$ be a non-trivial
\mbox{$t$-$(v,k,\lambda)$} design with $t \geq 2$. If $G \leq
\emph{Aut}(\mathcal{D})$ acts block-transitively on $\mathcal{D}$,
then $G$ acts point-transitively on $\mathcal{D}$.
\end{prop}

For a \mbox{$2$-$(v,k,1)$} design $\mathcal{D}$, it is elementary
that the point \mbox{$2$-transitivity} of \mbox{$G \leq
\mbox{Aut}(\mathcal{D})$} implies its flag-transitivity. For
\mbox{$2$-$(v,k,\lambda)$} designs, this implication remains true if
$r$ and $\lambda$ are relatively prime (see, for
instance,~\cite[Chap.\,2.3,\,Lemma\,8]{demb68}). However, for
\mbox{$t$-$(v,k,\lambda)$} designs with $t \geq 3$, it can be
deduced from Proposition~\ref{BlocksLemma2} that always the converse
holds (see~\cite{Buek1968} or~\cite[Lemma\,2]{Hu2001}):

\begin{prop}\label{flag2trs}
Let $\mathcal{D}=(X,\mathcal{B},I)$ be a non-trivial
\mbox{$t$-$(v,k,\lambda)$} design with \mbox{$t \geq 3$}. If
\mbox{$G \leq \emph{Aut}(\mathcal{D})$} acts flag-transitively on
$\mathcal{D}$, then $G$ acts point \mbox{$2$-transitively} on $\mathcal{D}$.
\end{prop}

Investigating highly symmetric \mbox{$t$-designs} for large values
of $t$, P.~Cameron and C.~Praeger~\cite[Thm.\,2.1]{CamPrae1993}
deduced from Theorem~\ref{RayCh} and Proposition~\ref{BlocksLemma2}
the following assertion:

\begin{prop}{\em (Cameron \& Praeger~1993).}\label{flag3hom}
Let $\mathcal{D}=(X,\mathcal{B},I)$ be a \mbox{$t$-$(v,k,\lambda)$}
design with $t\geq 2$. Then, the following holds:

\begin{enumerate}

\item[{\emph(a)}] If \mbox{$G \leq \emph{Aut}(\mathcal{D})$} acts block-transitively on $\mathcal{D}$,
then $G$ also acts point \mbox{$\lfloor t/2 \rfloor$-homogeneously} on $\mathcal{D}$.

\smallskip

\item[{\emph(b)}] If \mbox{$G \leq \emph{Aut}(\mathcal{D})$} acts flag-transitively on $\mathcal{D}$,
then $G$ also acts point \mbox{$\lfloor (t+1)/2 \rfloor$-homogeneously} on $\mathcal{D}$.

\end{enumerate}
\end{prop}

As for $t \geq 7$ the flag-transitivity, respectively for $t \geq 8$
the block-transitivity of \mbox{$G \leq \mbox{Aut} (\mathcal{D})$}
implies at least its point \mbox{$4$-homogeneity}, they obtained the
following restrictions as a consequence of the finite simple group
classification (\textit{cf.}~\cite[Thm.\,1.1]{CamPrae1993}):

\begin{thm}{\em (Cameron \& Praeger~1993).}
Let $\mathcal{D}=(X,\mathcal{B},I)$ be a \mbox{$t$-$(v,k,\lambda)$}
design. If \mbox{$G \leq \emph{Aut}(\mathcal{D})$} acts
block-transitively on $\mathcal{D}$ then $t \leq 7$, while if
\mbox{$G \leq \emph{Aut}(\mathcal{D})$} acts flag-transitively on
$\mathcal{D}$ then $t \leq 6$.
\end{thm}

Moreover, they formulated the following far-reaching conjecture
(\textit{cf.}~\cite[Conj.\,1.2]{CamPrae1993}):

\begin{conj}{\em (Cameron \& Praeger~1993).}
There are no non-trivial block-transitive \mbox{$6$-designs}.
\end{conj}

The author~\cite{Hu_mmics2008,Hu_cam2008} recently confirmed the
non-existence of block-transitive Steiner \mbox{$6$-designs}, modulo
two special cases that remain elusive.

\begin{thm}{\em (Huber~2010).}
Let $\mathcal{D}=(X,\mathcal{B},I)$ be a non-trivial Steiner
\mbox{$6$-design}. Then \mbox{$G \leq \emph{Aut}(\mathcal{D})$}
cannot act block-transitively on $\mathcal{D}$, except possibly when
$G= P \mathit{\Gamma} L(2,p^e)$ with $p=2$ or $3$ and $e$ is an odd
prime power.
\end{thm}

Previously, the author classified all flag-transitive Steiner
\mbox{$t$-designs} with $t>2$
(see~\cite{Hu2001,Hu2005,Hu_Habil2005,Hu2007,Hu2007a}
and~\cite{Hu2008} for a monograph). These results answered a series of
40-year-old problems and generalized theorems of
J.~Tits~\cite{Tits1964} and H.~L{\"u}neburg~\cite{Luene1965}. Earlier,
F.~Buekenhout, A.~Delandtsheer, J.~Doyen, P.~Kleidman, M.~Liebeck,
and J.~Saxl~\cite{Buek1990,Del2001,Kleid1990,Lieb1998,Saxl2002} had
essentially characterized all flag-transitive Steiner
\mbox{2-designs}. All these classification results rely on the
classification of the finite simple groups. As flag-transitivity
clearly implies block-transitivity, these results provide nice
examples of block-transitive Steiner designs. An encyclopedic survey
of further results, in particular on point-imprimitive
block-transitive \mbox{$t$-designs}, is given by~\cite{Del2007}.

\smallskip

%%% ----------------------------------------------------------------------

\section{Proof of the Main Theorem}\label{proof}

We prove in this section our Main Theorem stated in Section~\ref{intro}.
In order to investigate the existence problem of non-trivial
block-transitive Steiner \mbox{$7$-designs}, we can as a consequence
of Proposition~\ref{flag3hom}~(a) make use of the classification of
all finite \mbox{$3$-homogeneous} permutation groups, which itself
relies on the classification of all finite simple groups
(\textit{cf.}~\cite{Cam1981,Gor1982,Kant1972,Lieb1987,LivWag1965}).
We remark that the list given in~\cite[List\,2.2]{CamPrae1993} is slightly incomplete.

%%% ----------------------------------------------------------------------

\subsection{Finite $3$-homogeneous Permutation Groups}\label{list}

The list of groups is as follows: Let $G$ be a finite
\mbox{$3$-homogeneous} permutation group on a set $X$ with
$\left|X\right| \geq 4$. Then $G$ is either of

{\bf (A) Affine Type:} $G$ contains a regular normal subgroup $T$
which is elementary Abelian of order $v=2^d$. If we identify $G$
with a group of affine transformations
\[x \mapsto x^g+u\]
of $V=V(d,2)$, where $g \in G_0$ and $u \in V$, then particularly
one of the following occurs:

\begin{enumerate}

\smallskip

\item[(1)] $G \cong AGL(1,8)$, $A \mathit{\Gamma} L(1,8)$, or $A \mathit{\Gamma} L(1,32)$

\smallskip

\item[(2)] $G_0 \cong SL(d,2)$, $d \geq 2$

\smallskip

\item[(3)] $G_0 \cong A_7$, $v=2^4$

\end{enumerate}

or

\smallskip

{\bf (B) Almost Simple Type:} $G$ contains a simple normal subgroup
$N$, and \mbox{$N \leq G \leq \mbox{Aut}(N)$}. In particular, one of the
following holds, where $N$ and $v=|X|$ are given as follows:
\begin{enumerate}

\smallskip

\item[(1)] $A_v$, $v \geq 5$

\smallskip

\item[(2)] $PSL(2,q)$, $q>3$, $v=q+1$

\smallskip

\item[(3)] $M_v$, $v=11,12,22,23,24$ \hfill (Mathieu groups)

\smallskip

\item[(4)] $M_{11}$, $v=12$

\end{enumerate}

\smallskip

\noindent We note that if $q$ is odd, then $PSL(2,q)$ is
$3$-homogeneous for \mbox{$q \equiv 3$ (mod $4$)}, but not for
\mbox{$q \equiv 1$ (mod $4$)}, and hence not every group $G$ of
almost simple type satisfying (2) is $3$-homogeneous on $X$.

\begin{rem}\label{equa_t=5}
\emph{If \mbox{$G \leq \mbox{Aut}(\mathcal{D})$} acts
block-transitively on any Steiner \mbox{$t$-design} $\mathcal{D}$
with $t\geq 6$, then by Proposition~\ref{flag3hom}~(a), $G$ acts
point \mbox{$3$-homogeneously} and in particular point
\mbox{$2$-transitively} on $\mathcal{D}$. Applying
Corollary~\ref{Comb_t=5}~(b) yields the equation
\[b=\frac{{v \choose t}}{{k \choose
t}}=\frac{v(v-1) \left|G_{xy}\right|}{\left| G_B \right|},\] where
$x$ and $y$ are two distinct points in $X$ and $B$ is a block in
$\mathcal{B}$. We will see that this arithmetical condition in
combination with the combinatorial tools from Section~\ref{Prelim}
gives immediately strong results for some of the cases to be
examined.}
\end{rem}

%%% ---------------------------------------------------------------------

\subsection{Groups of Automorphisms of Affine Type}\label{affine
typ}

Using the notation as before, let us assume \emph{for the rest of
the section} that $\mathcal{D}=(X,\mathcal{B},I)$ is a non-trivial Steiner
\mbox{$7$-design} with \mbox{$G \leq \mbox{Aut}(\mathcal{D})$} acting
block-transitively on $\mathcal{D}$. Clearly, we may assume that $k>7$ as we
do not consider trivial designs. We will examine in this subsection
those cases where $G$ is of affine type.

\medskip

\emph{Case} (1): $G \cong AGL(1,8)$, $A \mathit{\Gamma} L(1,8)$, or
$A \mathit{\Gamma} L(1,32)$.

\smallskip

As $k>7$, the case $v=8$ is not possible. For $v=32$, we have
$\left| G \right| =5v(v-1)$ and $k \leq 11$ by
Corollary~\ref{Cameron_t=5}. The few possibilities can easily be
ruled out using Corollary~\ref{Comb_t=5} together with
Remark~\ref{equa_t=5}.

\medskip

\emph{Case} (2): $G_0 \cong SL(d,2)$, $d \geq 2$.

\smallskip

Let $e_i$ denote the $i$-th standard basis vector of the vector
space $V=V(d,2)$, and $\text{\footnotesize{$\langle$}} e_i
\text{\footnotesize{$\rangle$}}$ the \mbox{$1$-dimensional} vector
subspace spanned by $e_i$. We will prove by contradiction that
\mbox{$G \leq \mbox{Aut}(\mathcal{D})$} cannot act block-transitively on any
non-trivial Steiner \mbox{$7$-design} $\mathcal{D}$.

We may assume that $v=2^d > k > 7$. We remark that clearly any seven
distinct points are non-coplanar in $AG(d,2)$ and hence generate an
affine subspace of dimension at least $3$. Let $\mathcal{E}=
\text{\footnotesize{$\langle$}} e_1,e_2,e_3
\text{\footnotesize{$\rangle$}}$ denote the \mbox{$3$-dimensional}
vector subspace spanned by $e_1,e_2,e_3$. Then
$SL(d,2)_\mathcal{E}$, and therefore also $G_{0,\mathcal{E}}$, acts
point-transitively on \mbox{$V \setminus \mathcal{E}$}. If the
unique block $B \in \mathcal{B}$ which is incident with the \mbox{$7$-subset}
\mbox{$\{0,e_1,e_2,e_3,e_1+e_2,e_2+e_3, e_1+e_3\}$} contains some
point outside $\mathcal{E}$, then it would already contain all
points of \mbox{$V \setminus \mathcal{E}$}. But then \mbox{$k \geq
2^d-8+7=2^d-1$}, a contradiction to Corollary~\ref{Cameron_t=5}.
Hence, $B$ lies completely in $\mathcal{E}$. The block-transitivity
of $G$ now implies that each block must be contained in a
\mbox{$3$-dimensional} affine subspace. On the other hand, any
seven distinct points that do not lie in a \mbox{$3$-dimensional} affine
subspace must also be incident with a unique block by the definition
of Steiner designs, a contradiction.

\medskip

\emph{Case} (3): $G_0 \cong A_7$, $v=2^4$.

\smallskip

As $v=2^4$, we have $k \leq 9$ by Corollary~\ref{Cameron_t=5}. But,
Corollary~\ref{Comb_t=5}~(c) obviously eliminates the cases when
$k=8$ or $9$.

%%% ---------------------------------------------------------------------

\subsection{\mbox{Groups of Automorphisms of Almost Simple Type}}\label{almost simple
type}

When $G$ is of almost simple type, then the Cases (B) (1), (3) and (4) of Section~\ref{list}
cannot occur as elementarily proved in~\cite[Sect.\,2,\;\mbox{mainly}\;Prop.\,2.4]{CamPrae1993}.
Hence, we only have to consider

\medskip

\emph{Case} (2): $N=PSL(2,q)$, $v=q+1$, $q=p^e >3$.

\smallskip

Here $\mbox{Aut}(N)= P \mathit{\Gamma} L (2,q)$, and $\left| G \right| =
(q+1)q \frac{(q-1)}{n}a$ with $n=(2,q-1)$ and $a \mid ne$. We may
assume that $q \geq 8$. We will show that \mbox{$G \leq \mbox{Aut}(\mathcal{D})$}
cannot act block-transitively on any non-trivial Steiner
\mbox{$7$-design} $\mathcal{D}$.

From Remark~\ref{equa_t=5}, we obtain
\begin{equation}\label{Eq-0}
(q-2)(q-3)(q-4)(q-5) \left| G_{B} \right| \cdot n =
k(k-1)(k-2)(k-3)(k-4)(k-5)(k-6) \cdot a.
\end{equation}
Due to Proposition~\ref{Cam}~(b), we have
\begin{equation}\label{Eq-A}
q-5 \geq (k-5)(k-6).
\end{equation}
Hence, from Equation~(\ref{Eq-0}) follows
\begin{equation}\label{Eq-B}
(q-2)(q-3)(q-4)\left|G_{B} \right| \cdot n \leq
k(k-1)(k-2)(k-3)(k-4) \cdot a.
\end{equation}
If we assume that $k \geq 27$, then obviously
\[k(k-1)(k-2)(k-3) < 2[(k-5)(k-6)]^2,\]
and hence
\[(q-2)(q-3)(q-4)\left|G_{B} \right| \cdot n< 2(q-5)^2 (k-4)\cdot a \leq 2(q-5)^2 \bigl\lfloor \sqrt{q+1} + \frac{3}{2}\bigr\rfloor  \cdot a\]
in view of Inequality~(\ref{Eq-A}) and Corollary~\ref{Cameron_t=5}.
Taking into account that always $a \leq$ log$_2 q$, it follows
immediately that \mbox{$\left|G_{B} \right| \cdot n=1$} for each
possible value of $q \neq 32$. Hence, in particular $q$ must be
even. But then the right hand side of Equation~(\ref{Eq-0}) is
always divisible by $16$ but never the left hand side, a
contradiction. The case $q=32$ as well as the few remaining
possibilities for $k<27$ can easily be ruled out by hand using
Equation~(\ref{Eq-0}) and Inequality~(\ref{Eq-A}).

\smallskip

This completes the proof of the Main Theorem.\qed

\acknowledgements

The author gratefully acknowledges support by the Deutsche
Forschungsgemeinschaft (DFG) via a Heisenberg grant (Hu954/4).

\end{document}